\documentclass[12pt,reqno]{amsart}
\usepackage[margin=1in]{geometry}
\usepackage{amsmath,amssymb,amsthm}
\usepackage[colorlinks=true,linkcolor=blue,citecolor=blue,urlcolor=blue]{hyperref}

\usepackage{etoolbox}

\makeatletter
\patchcmd{\@setauthors}
  {\MakeUppercase{\authors}}
  {\authors}
  {}{}
\makeatother

\makeatletter
\patchcmd{\@settitle}
  {\uppercasenonmath\@title}
  {}{}{}
\makeatother

\newtheorem{theorem}{Theorem}[section]
\newtheorem{lemma}[theorem]{Lemma}
\newtheorem{proposition}[theorem]{Proposition}
\newtheorem{corollary}[theorem]{Corollary}
\theoremstyle{remark}
\newtheorem{example}[theorem]{Example}
\newtheorem{remark}[theorem]{Remark}

\numberwithin{equation}{section}

\newcommand{\one}{\mathbf 1}
\newcommand{\ee}{\mathbf e}
\newcommand{\rhL}{\rho_L}
\DeclareMathOperator{\diag}{diag}
\DeclareMathOperator{\Span}{span}

\allowdisplaybreaks[1]

\begin{document}

\title{Graphs attaining Haemers' toughness bound}

\author{
Yifan Li\textsuperscript{*}
\qquad
Cheng Yeaw Ku\textsuperscript{$\dagger$}
}

\thanks{
\textsuperscript{*}
 School of Mathematics and Statistics,
Fujian Normal University,
Fuzhou 350117, P. R. China.
Email: \texttt{llyfliyifan@163.com}
}

\thanks{
\textsuperscript{$\dagger$}
 Division of Mathematical Sciences,
Nanyang Technological University,
Singapore 637371, Singapore.
Email: \texttt{cyku@ntu.edu.sg}
}

\date{}

\begin{abstract}
Haemers \cite{Haemers2020} conjectured that every connected noncomplete graph $G$ satisfies
$t(G)\geqslant\mu_2/(\mu_n-\delta)$, where $t(G)$ is the toughness of
$G$, $\mu_2$ and $\mu_n$ are its algebraic connectivity and Laplacian
spectral radius, and $\delta$ is its minimum degree. Greaves and Zhu \cite{GreavesZhu2026}
recently proved this conjecture, leaving the equality case open. We characterize all graphs attaining
equality: a graph attains equality if and only if its complement has a
complete component whose order is at least the Laplacian spectral radius
of every other component. Equivalently, its algebraic connectivity and
independence number sum to its order. We also determine all vertex cuts
attaining the toughness, and show that for regular graphs, the equality graphs are
exactly the complete multipartite graphs with parts of equal size. 
\end{abstract}

\maketitle
\pagestyle{plain}
\section{Introduction}\label{sec:introduction}

\subsection{Toughness and Laplacian bounds}\label{subsec:background}

All graphs in this paper are finite, simple, and undirected. For a graph
$G=(V,E)$ and a subset $U\subseteq V$, let $G[U]$ denote the subgraph
induced by $U$. For $S\subseteq V$, write $G-S=G[V\setminus S]$, and let
$c(G-S)$ be the number of connected components of $G-S$. A set $S$ is a
\emph{vertex cut} if $c(G-S)\geqslant2$, and the \emph{vertex
connectivity} $\kappa(G)$ of a noncomplete graph $G$ is the minimum cardinality among all vertex cuts of \(G\). The \emph{toughness} of a connected noncomplete graph
$G$ is
\begin{equation}\label{eq:toughness-definition}
 t(G)=\min\left\{\frac{|S|}{c(G-S)}:
 S\subseteq V,\ c(G-S)\geqslant2\right\},
\end{equation}
and a complete graph has infinite toughness by convention. A vertex cut
attaining the minimum in \eqref{eq:toughness-definition} is called a
\emph{toughness cut}.

Chv\'atal~\cite{Chvatal1973} introduced toughness in his study of
Hamiltonian graphs. While vertex connectivity is determined by the minimum size of a vertex cut, toughness counts the number of components resulting from vertex deletion.
 If $C$ is a Hamiltonian cycle of $G$ and $S$ is a
nonempty vertex set, then $C-S$ is a spanning subgraph of $G-S$ with at
most $|S|$ components. Hence every Hamiltonian graph is
\emph{$1$-tough}, that is, $t(G)\geqslant1$. Toughness has since been
studied in connection with Hamiltonicity, factors, matchings, and
spanning trees; see the survey of Bauer, Broersma, and
Schmeichel~\cite{BauerBroersmaSchmeichel2006}.

Spectral bounds control toughness without examining individual vertex
cuts. Alon~\cite{Alon1995} and Brouwer~\cite{Brouwer1995} initiated this
approach for regular graphs. For a connected $k$-regular graph $G$, let
$\lambda$ be the largest absolute value of an adjacency eigenvalue other
than $k$. Brouwer proved that $t(G)>k/\lambda-2$ and conjectured that
$t(G)\geqslant k/\lambda-1$; Gu~\cite{Gu2021} proved this conjecture.
Cioab\u{a} and Wong~\cite{CioabaWong2014} and Cioab\u{a} and
Gu~\cite{CioabaGu2016} further developed the relations among the
adjacency spectrum, toughness, and other connectivity parameters.

For irregular graphs, the Laplacian matrix $L(G)=D(G)-A(G)$ provides a natural framework for spectral analysis,
 where $A(G)$ is the adjacency matrix
of $G$ and $D(G)$ is the diagonal matrix of its vertex degrees. For a
graph of order $n$, the Laplacian eigenvalues of \(G\) are denoted by
\[
 0=\mu_1(G)\leqslant\mu_2(G)\leqslant\cdots\leqslant\mu_n(G),
\]
and we omit the argument $G$ when the graph is clear from the context.
The eigenvalue $\mu_2$ is the \emph{algebraic connectivity} introduced by
Fiedler~\cite{Fiedler1973}, and $\delta=\delta(G)$ denotes the minimum
degree. Gu and Haemers~\cite{GuHaemers2022} proved that every connected
noncomplete graph $G$ of order $n$ satisfies
\begin{align}
 t(G)&\geqslant\frac{\mu_n\mu_2}{n(\mu_n-\delta)},
       \label{eq:gu-haemers-first}\\
 t(G)&\geqslant\frac{\mu_2}{\mu_n-\mu_2},
       \label{eq:gu-haemers-second}
\end{align}
and Huang, Das, and Zhu~\cite{HuangDasZhu2022} characterized the graphs
attaining equality in each of these bounds.

Haemers \cite{Haemers2020} (see also \cite[Conjecture~2.1]{GuHaemers2022}) conjectured the stronger
inequality
\begin{equation}\label{eq:haemers-bound}
 t(G)\geqslant\frac{\mu_2}{\mu_n-\delta}.
\end{equation}
Since $\mu_2\leqslant\delta$ and $\mu_n\leqslant n$
(Lemma~\ref{lem:basic-bounds}), inequality \eqref{eq:haemers-bound}
implies both \eqref{eq:gu-haemers-first} and
\eqref{eq:gu-haemers-second}. Greaves and Zhu~\cite{GreavesZhu2026}
proved \eqref{eq:haemers-bound} and leaving the characterization of
the graphs attaining equality open. As noted by Gu and Haemers \cite{GuHaemers2022}, every complete multipartite graph with a part of size at least two attains equality in \eqref{eq:haemers-bound}. Greaves and Zhu \cite{GreavesZhu2026} further exhibited non-complete multipartite equality graphs obtained by adding edges within certain parts other than a largest part. These constructions are subsumed by Example~\ref{ex:multipartite} and
Proposition~\ref{prop:edge-addition} below. In this paper, we characterize all
graphs attaining equality.

\subsection{Main results}\label{subsec:results}

Let $\overline G$ denote the complement of $G$,
$K_m$ the complete graph on $m$ vertices,
and $F_1\sqcup F_2$ the disjoint union of $F_1$ and $F_2$. The \emph{join} $F_1\vee F_2$ is defined as
the disjoint union $F_1\sqcup F_2$
together with all edges joining
$V(F_1)$ and $V(F_2)$.
The \emph{Laplacian spectral radius} $\rhL(F)$
of a nonempty graph $F$ is the largest
eigenvalue of $L(F)$, thus
$\rhL(F)=0$ if and only if $F$ has no edges. Finally, $\alpha(G)$ denotes
the independence number of $G$, the maximum size of a set of pairwise
nonadjacent vertices.

\begin{theorem}\label{thm:main}
Let $G$ be a connected noncomplete graph of order $n$ with minimum
degree $\delta$ and Laplacian eigenvalues
$0=\mu_1<\mu_2\leqslant\cdots\leqslant\mu_n$.
The following statements are equivalent.
\begin{enumerate}
 \item[\textup{(i)}] $t(G)=\mu_2/(\mu_n-\delta)$.
 \item[\textup{(ii)}] $\mu_2+\alpha(G)=n$.
 \item[\textup{(iii)}] There exist an integer $m\geqslant2$ and a graph
 $F$ with at least one vertex such that
 \[
   \overline G\cong K_m\sqcup F
   \qquad\text{and}\qquad
   \rhL(F)\leqslant m.
 \]
\end{enumerate}
If these statements hold, then the integer $m$ in \textup{(iii)} equals
$\alpha(G)$, and
\begin{equation}\label{eq:equality-parameters}
 \mu_2=\delta=n-m,\qquad
 \mu_n=n,\qquad
 t(G)=\frac{n-m}{m}.
\end{equation}
\end{theorem}

We call a connected noncomplete graph an \emph{equality graph} if it
attains equality in \eqref{eq:haemers-bound}.

The three conditions in Theorem~\ref{thm:main} describe equality from
complementary points of view. Condition~\textup{(ii)} is the equality
case of the chain
\begin{equation}\label{eq:alpha-chain}
 \mu_2\leqslant\delta\leqslant n-\alpha(G),
\end{equation}
which holds for every connected noncomplete graph: the first inequality
is part of Lemma~\ref{lem:basic-bounds}, and the second holds because a
vertex in a maximum independent set has all its neighbors outside that
set. Condition~\textup{(iii)} is structural. It states that
$G\cong\overline K_m\vee H$ with $H=\overline F$: the graph $G$ contains
an independent set of size $m$ whose vertices are adjacent to all other
vertices, and the remaining part $H$ is constrained only through the
Laplacian spectral radius of its complement. In particular, the
components of $F$ need not be complete, and their orders are not bounded
in terms of $m$ (Examples~\ref{ex:cycle-five} and~\ref{ex:paths}).

Regularity removes this freedom. For positive integers
$n_1,\ldots,n_r$, the complete multipartite graph $K_{n_1,\ldots,n_r}$
has independent parts of sizes $n_1,\ldots,n_r$ and all edges between
distinct parts.

\begin{corollary}\label{cor:regular-introduction}
A connected noncomplete regular graph is an equality graph if and only if
it is a complete multipartite graph whose parts have equal size, that
is, $K_{m,\ldots,m}$ with $r$ parts for some $r,m\geqslant2$. Its
toughness is then $r-1$.
\end{corollary}

The proof of Theorem~\ref{thm:main} also determines all toughness cuts
of an equality graph: they are exactly the complements of its maximum
independent sets, or equivalently, the complements of the components of
$\overline G$ isomorphic to $K_{\alpha(G)}$ (Theorem~\ref{thm:all-cuts}).
Consequently, every equality graph satisfies
$\mu_2(G)=\kappa(G)=\delta(G)$, and each of its toughness cuts is a
minimum vertex cut (Corollary~\ref{cor:minimum-cuts}).
Fiedler~\cite{Fiedler1973} proved that $\mu_2(G)\leqslant\kappa(G)$ for
every noncomplete graph, and Kirkland, Molitierno, Neumann, and
Shader~\cite{KirklandMolitiernoNeumannShader2002} characterized the
graphs attaining equality; every equality graph for
\eqref{eq:haemers-bound} therefore belongs to the class they describe.
A minimum vertex cut of an equality graph need not be a toughness cut
(Example~\ref{ex:minimum-not-toughness}).

Finally, \eqref{eq:haemers-bound} is the strongest of the three bounds,
so equality in \eqref{eq:gu-haemers-first} or in
\eqref{eq:gu-haemers-second} implies equality in
\eqref{eq:haemers-bound}. Theorem~\ref{thm:main} yields the converse:
every equality graph satisfies $\mu_2=\delta$ and $\mu_n=n$, and hence
attains equality in all three bounds (Corollary~\ref{cor:three-bounds}).
Thus the three bounds have the same equality cases, and
Theorem~\ref{thm:main}\textup{(iii)} is an equivalent form of the
equality family found by Huang, Das, and
Zhu~\cite[Theorem~1.3]{HuangDasZhu2022}; see
Section~\ref{subsec:join}.

\subsection{Outline of the proof}\label{subsec:outline}

Fix a vertex cut $S$, let $H_1,\ldots,H_c$ be the components of $G-S$,
and let $r_i$ be the average number of neighbors in $S$ of a vertex of
$H_i$. Compressing $L(G)$ onto the span of the characteristic vectors of
$V(H_1),\ldots,V(H_c)$ and $S$ yields two inequalities. The first, which
uses $\mu_n$, gives $|S|(\mu_n-\delta)\geqslant r_1+\cdots+r_c$
(Lemma~\ref{lem:upper-cut}). The second, which uses $\mu_2$, shows that
at most one $r_i$ is smaller than $\mu_2$ (Lemma~\ref{lem:lower-cut}).
The key step is Proposition~\ref{prop:exceptional-strictness}: if some
$r_i$ is smaller than $\mu_2$, then the ratio $|S|/c$ is strictly larger
than $\mu_2/(\mu_n-\delta)$. Therefore a cut attaining equality has
$r_i=\mu_2$ for every $i$, and the equality cases of
Lemma~\ref{lem:upper-cut} force every $H_i$ to be a single vertex of
degree $\delta=\mu_2$ (Proposition~\ref{prop:cut-bound-equality}).

Now let $S$ be a toughness cut of an equality graph. The vertices outside
$S$ form an independent set, and each of them has degree $\mu_2$. A
Rayleigh quotient argument shows that they have a common neighborhood
(Lemma~\ref{lem:common-neighborhood}), and the minimality of $S$ forces
this neighborhood to be all of $S$
(Lemma~\ref{lem:equality-cut-structure}). This gives the decomposition
of $\overline G$ in Theorem~\ref{thm:main}\textup{(iii)}. Since the
argument applies to every toughness cut, it determines all of them.

The paper is organized as follows. In Section~\ref{sec:preliminaries}
we collect the Laplacian facts used in the proofs.
In Section~\ref{sec:cuts} we prove the vertex-cut inequality together with its equality case.
In Section~\ref{sec:characterization} we prove Theorem~\ref{thm:main} and
compares the equality cases of the three bounds.
In Section~\ref{sec:consequences} we treat regular graphs and toughness cuts.
In Section~\ref{sec:examples} we present constructions, the equality graphs
with independence number two or three, and a recognition criterion.

\section{Laplacian preliminaries}\label{sec:preliminaries}

For a vertex $v$ of a graph $F$, we write $d_F(v)$ for its degree,
$N_F(v)$ for its neighborhood, and $N_F[v]=N_F(v)\cup\{v\}$ for its
closed neighborhood. For disjoint vertex sets $U$ and $W$, let
$e_F(U,W)$ be the number of edges with one end in $U$ and the other in
$W$. The maximum degree of $F$ is $\Delta(F)$. We omit the subscript $F$
when the graph is clear from the context.

Vectors are real column vectors indexed by the vertices of the graph
under consideration. We write $\one$ for the all-ones vector,
$J=\one\one^{\mathsf T}$ for the all-ones matrix, $\ee_v$ for the
coordinate vector of a vertex $v$, and $\chi(U)$ for the characteristic
vector of a vertex set $U$. The identity matrix is $I$, with its order
determined by the context. For every vector $x$,
\begin{equation}\label{eq:laplacian-energy}
 x^{\mathsf T}L(F)x=\sum_{uv\in E(F)}(x_u-x_v)^2.
\end{equation}
Hence $L(F)$ is positive semidefinite, and its kernel consists of the
vectors that are constant on each connected component of $F$. For a
graph $G$ of order $n$,
\begin{equation}\label{eq:variational}
 \mu_2(G)=\min_{\substack{x\neq0\\x\perp\one}}
       \frac{x^{\mathsf T}L(G)x}{x^{\mathsf T}x},
 \qquad
 \mu_n(G)=\max_{x\neq0}
       \frac{x^{\mathsf T}L(G)x}{x^{\mathsf T}x}.
\end{equation}
Expanding in an orthonormal eigenbasis of $L(G)$ shows that every nonzero
vector $x\perp\one$ attaining the minimum in \eqref{eq:variational} is a
$\mu_2(G)$-eigenvector.

These facts are standard; see Fiedler~\cite{Fiedler1973} and Brouwer and
Haemers~\cite{BrouwerHaemers2012}. We include short proofs of the
consequences below because their equality cases are used later.

\begin{lemma}\label{lem:complement}
Let $F$ be a graph of order $q\geqslant2$. Then
\begin{equation}\label{eq:complement-identity}
 L(F)+L(\overline F)=qI-J
 \qquad\text{and}\qquad
 \mu_2(F)=q-\rhL(\overline F).
\end{equation}
Moreover, $\rhL(F)\leqslant q$, with equality if and only if
$\overline F$ is disconnected.
\end{lemma}

\begin{proof}
The matrix identity follows by comparing diagonal and off-diagonal
entries. Both Laplacian matrices annihilate $\one$, and on $\one^\perp$
their sum acts as $qI$. Hence the eigenvalues of $L(\overline F)$ on
$\one^\perp$ are $q-\mu_i(F)$ for $2\leqslant i\leqslant q$, and the
largest of them is $q-\mu_2(F)\geqslant0$. The remaining eigenvalue of
$L(\overline F)$, on $\Span\{\one\}$, is zero, so
$\rhL(\overline F)=q-\mu_2(F)$.

Since $L(\overline F)$ is positive semidefinite, every eigenvalue of
$L(F)$ on $\one^\perp$ is at most $q$, and equality occurs if and only if
$L(\overline F)$ has a nonzero null vector orthogonal to $\one$. By
\eqref{eq:laplacian-energy}, this happens if and only if $\overline F$
has at least two connected components.
\end{proof}

\begin{lemma}\label{lem:degree-radius}
Let $F$ be a graph with at least one edge. Then
\begin{equation}\label{eq:degree-radius}
 \Delta(F)+1\leqslant\rhL(F)\leqslant2\Delta(F).
\end{equation}
More generally, $\rhL(F)\geqslant d_F(v)+1$ for every vertex $v$ of
positive degree, and this inequality is strict if some edge joins
$N_F(v)$ to $V(F)\setminus N_F[v]$.
\end{lemma}

\begin{proof}
Fix a vertex $v$ of positive degree $d=d_F(v)$. Define $x$ by $x_v=d$,
$x_u=-1$ for $u\in N_F(v)$, and $x_w=0$ for all other vertices $w$.
Then $x^{\mathsf T}x=d(d+1)$, and \eqref{eq:laplacian-energy} gives
\[
 x^{\mathsf T}L(F)x
 =d(d+1)^2+e_F\bigl(N_F(v),V(F)\setminus N_F[v]\bigr).
\]
The lower bounds and the strictness assertion follow from the maximum
Rayleigh quotient; for the first inequality in \eqref{eq:degree-radius},
choose $v$ of maximum degree.

For the upper bound, \eqref{eq:laplacian-energy} gives, for every $x$,
\[
 x^{\mathsf T}L(F)x
 \leqslant2\sum_{uv\in E(F)}(x_u^2+x_v^2)
 =2\sum_{v\in V(F)}d_F(v)x_v^2
 \leqslant2\Delta(F)\,x^{\mathsf T}x.
\]
Taking the maximum Rayleigh quotient completes the proof.
\end{proof}

\begin{lemma}\label{lem:basic-bounds}
Let $G$ be a connected noncomplete graph of order $n$ with minimum
degree $\delta$. Then
\begin{equation}\label{eq:basic-bounds}
 0<\mu_2\leqslant\delta<\mu_n-1\leqslant n-1.
\end{equation}
Moreover, $\mu_2\leqslant|S|$ for every vertex cut $S$ of $G$; that is,
$\mu_2\leqslant\kappa(G)$.
\end{lemma}

\begin{proof}
Since $G$ is connected, the kernel of $L(G)$ is spanned by $\one$, and
therefore $\mu_2>0$.

Let $S$ be a vertex cut, and let $U$ and $W$ be the vertex sets of two
distinct components of $G-S$. Define $x$ to be $|W|$ on $U$, $-|U|$ on
$W$, and $0$ elsewhere. Then $x\neq0$ and $x\perp\one$. By
\eqref{eq:laplacian-energy}, only edges between $U\cup W$ and $S$
contribute to $x^{\mathsf T}L(G)x$, and every vertex of $U\cup W$ has at
most $|S|$ neighbors in $S$. Hence
$x^{\mathsf T}L(G)x\leqslant|S|\,x^{\mathsf T}x$, and
\eqref{eq:variational} gives $\mu_2\leqslant|S|$.

Let $v$ be a vertex of degree $\delta$. Since $G$ is noncomplete,
$\delta\leqslant n-2$. Deleting $N(v)$ therefore isolates $v$ and leaves
at least one further vertex, so $N(v)$ is a vertex cut and
$\mu_2\leqslant\delta$. Since $G$ is connected and
$V(G)\setminus N[v]$ is nonempty, some edge joins $N(v)$ to
$V(G)\setminus N[v]$, and Lemma~\ref{lem:degree-radius} gives
$\mu_n>\delta+1$. Finally, $\mu_n\leqslant n$ by
Lemma~\ref{lem:complement}.
\end{proof}

\begin{lemma}\label{lem:common-neighborhood}
Let $G$ be a connected graph, and let $u$ and $v$ be distinct
nonadjacent vertices with $d_G(u)=d_G(v)=\mu_2(G)$. Then
$N_G(u)=N_G(v)$.
\end{lemma}

\begin{proof}
The vector $x=\ee_u-\ee_v$ is orthogonal to $\one$. Since $u$ and $v$
are not adjacent,
\[
 \frac{x^{\mathsf T}L(G)x}{x^{\mathsf T}x}
 =\frac{d_G(u)+d_G(v)}2=\mu_2(G).
\]
Thus $x$ attains the minimum in \eqref{eq:variational}, and
$L(G)x=\mu_2(G)x$. At a vertex $w\notin\{u,v\}$, this equation reads
$-A(G)_{wu}+A(G)_{wv}=0$. Hence every such $w$ is adjacent to $u$ if
and only if it is adjacent to $v$.
\end{proof}

\section{Vertex cuts and the equality case}\label{sec:cuts}

Throughout this section, $G$ is a connected noncomplete graph of order
$n$ with minimum degree $\delta$ and Laplacian eigenvalues
$0=\mu_1<\mu_2\leqslant\cdots\leqslant\mu_n$. Fix a vertex cut $S$, let
$H_1,\ldots,H_c$ be the components of $G-S$, where
$c=c(G-S)\geqslant2$, and set
\begin{equation}\label{eq:cut-notation}
 h_i=|V(H_i)|,
 \qquad
 r_i=\frac{e_G(V(H_i),S)}{h_i}
 \qquad (1\leqslant i\leqslant c).
\end{equation}
Thus $r_i$ is the average number of neighbors in $S$ of a vertex of
$H_i$. Since $G$ is connected, $r_i>0$. A vertex of $H_i$ has at most
$h_i-1$ neighbors outside $S$, so
\begin{equation}\label{eq:boundary-degree}
 r_i\geqslant\delta-h_i+1.
\end{equation}

Lemmas~\ref{lem:upper-cut} and~\ref{lem:lower-cut} are the compression
inequalities underlying the proof of Greaves and
Zhu~\cite[Lemmas~2.3 and~2.4]{GreavesZhu2026}; we derive them in a form
that displays their equality cases. The new ingredient is
Proposition~\ref{prop:exceptional-strictness}, which shows that a
component with $r_i<\mu_2$ forces strict inequality.

Let $T$ be the $n\times(c+1)$ matrix whose $i$th column is
$\chi(V(H_i))/\sqrt{h_i}$ for $1\leqslant i\leqslant c$ and whose last
column is $\chi(S)/\sqrt{|S|}$. These columns are orthonormal, and a direct
computation shows that the compression $Q=T^{\mathsf T}L(G)T$ is
\begin{equation}\label{eq:quotient}
 Q=
 \begin{pmatrix}
 r_1 & & 0 & -r_1\sqrt{h_1/|S|}\\
     &\ddots& &\vdots\\
 0   & & r_c & -r_c\sqrt{h_c/|S|}\\
 -r_1\sqrt{h_1/|S|}&\cdots&-r_c\sqrt{h_c/|S|}
      &\displaystyle\frac{1}{|S|}\sum_{i=1}^c h_ir_i
 \end{pmatrix}.
\end{equation}
The matrices $\mu_nI-L(G)$ and $L(G)-\mu_2I+(\mu_2/n)J$ are positive
semidefinite; the second vanishes on $\one$ and agrees with
$L(G)-\mu_2I$ on $\one^\perp$. Hence their compressions by $T$ are
positive semidefinite as well.

\begin{lemma}\label{lem:upper-cut}
In the notation of \eqref{eq:cut-notation}, $r_i<\mu_n$ for every $i$,
and
\begin{equation}\label{eq:schur-bound}
 |S|\geqslant\sum_{i=1}^c\frac{h_ir_i}{\mu_n-r_i}.
\end{equation}
Moreover, for every $i$,
\begin{equation}\label{eq:denominator-bound}
 \mu_n-r_i
 \leqslant\mu_n-\delta+h_i-1
 \leqslant(\mu_n-\delta)h_i,
\end{equation}
where the second inequality is strict if and only if $h_i\geqslant2$. Consequently,
\begin{equation}\label{eq:sum-bound}
 |S|(\mu_n-\delta)\geqslant\sum_{i=1}^c r_i.
\end{equation}
\end{lemma}

\begin{proof}
The matrix $\mu_nI-Q$, which is the compression of $\mu_n I - L(G)$ by $T$, is positive semidefinite. For $1 \le i \le c$, its $i$th diagonal entry
is $\mu_n-r_i$, and its $(i,c+1)$ entry $r_i\sqrt{h_i/|S|}$ is nonzero.
In a positive semidefinite matrix, a zero diagonal entry forces the
corresponding row to vanish, as the principal $2\times2$ submatrices
show. Hence $\mu_n-r_i>0$.

The leading principal $c\times c$ submatrix of $\mu_nI-Q$ is therefore
diagonal and positive definite, and its Schur complement in
$\mu_nI-Q$ is nonnegative:
\[
 0
 \leqslant\mu_n-\frac1{|S|}\sum_{i=1}^c h_ir_i
       -\frac1{|S|}\sum_{i=1}^c\frac{h_ir_i^2}{\mu_n-r_i}
 =\mu_n\left(1-\frac1{|S|}
               \sum_{i=1}^c\frac{h_ir_i}{\mu_n-r_i}\right).
\]
This proves \eqref{eq:schur-bound}.

The first inequality in \eqref{eq:denominator-bound} is
\eqref{eq:boundary-degree}. For the second,
\[
 (\mu_n-\delta)h_i-(\mu_n-\delta+h_i-1)=(h_i-1)(\mu_n-\delta-1),
\]
which is nonnegative, and positive if and only if $h_i\geqslant2$, because
$\mu_n-\delta>1$ by Lemma~\ref{lem:basic-bounds}. Finally,
\eqref{eq:denominator-bound} gives
$h_ir_i/(\mu_n-r_i)\geqslant r_i/(\mu_n-\delta)$ for every $i$, and
substituting these inequalities into \eqref{eq:schur-bound} yields
\eqref{eq:sum-bound}.
\end{proof}

\begin{lemma}\label{lem:lower-cut}
In the notation of \eqref{eq:cut-notation}, any two distinct indices
$i$ and $j$ satisfy
\begin{equation}\label{eq:pair-bound}
 (r_i-\mu_2)h_j+(r_j-\mu_2)h_i\geqslant0.
\end{equation}
In particular, at most one of $r_1,\ldots,r_c$ is smaller than $\mu_2$.
If $r_1<\mu_2$, then $r_j>\mu_2$ for every $2\leqslant j\leqslant c$,
and
\begin{equation}\label{eq:exceptional-bound}
 n+\mu_2\sum_{j=2}^c\frac{h_j}{r_j-\mu_2}
 \leqslant\frac{\mu_2h_1}{\mu_2-r_1}.
\end{equation}
\end{lemma}

\begin{proof}
The leading principal $c\times c$ submatrix of the compression of
$L(G)-\mu_2I+(\mu_2/n)J$ is
\[
 X=\diag(r_1-\mu_2,\ldots,r_c-\mu_2)
       +\frac{\mu_2}{n}uu^{\mathsf T},
 \qquad
 u=(\sqrt{h_1},\ldots,\sqrt{h_c})^{\mathsf T},
\]
and it is positive semidefinite. For distinct $i$ and $j$, evaluate its
quadratic form at the vector with $i$th coordinate $\sqrt{h_j}$, $j$th
coordinate $-\sqrt{h_i}$, and all other coordinates zero. This vector is
orthogonal to $u$, so the rank-one term vanishes, and
\eqref{eq:pair-bound} follows. The assertions about the signs of
$r_i-\mu_2$ follow immediately.

Now suppose that $r_1<\mu_2$. Then every diagonal entry $r_i-\mu_2$ is
nonzero, and exactly one of them is negative. By the determinant formula
for a rank-one perturbation of a diagonal matrix,
\[
 \det X
 =\prod_{i=1}^c(r_i-\mu_2)
   \left(1+\frac{\mu_2}{n}
          \sum_{i=1}^c\frac{h_i}{r_i-\mu_2}\right).
\]
Since $\det X\geqslant0$ and the product is negative, the factor in
parentheses is nonpositive. Separating the term $i=1$ and multiplying by
$n$ gives \eqref{eq:exceptional-bound}.
\end{proof}

\begin{proposition}\label{prop:exceptional-strictness}
In the notation of \eqref{eq:cut-notation}, if $r_i<\mu_2$ for some $i$,
then
\begin{equation}\label{eq:strict-cut}
 |S|(\mu_n-\delta)>c\mu_2.
\end{equation}
\end{proposition}

\begin{proof}
Renumber the components so that $r_1<\mu_2$, and put
\[
 a=\mu_2-r_1>0,
 \qquad
 b_j=r_j-\mu_2\quad(2\leqslant j\leqslant c).
\]
By Lemma~\ref{lem:lower-cut}, $b_j>0$ for every $j$. Suppose, for a
contradiction, that $|S|(\mu_n-\delta)\leqslant c\mu_2$. Combined with
\eqref{eq:sum-bound}, this gives $\sum_{i=1}^c r_i\leqslant c\mu_2$,
that is,
\begin{equation}\label{eq:loss-budget}
 a\geqslant\sum_{j=2}^c b_j.
\end{equation}
In particular, $a\geqslant b_j$ for every $j\geqslant2$.

Substitute $n=|S|+\sum_{i=1}^c h_i$ into \eqref{eq:exceptional-bound}
and multiply by $a$. Since $\mu_2h_1-ah_1=r_1h_1$ and
$1+\mu_2/b_j=r_j/b_j$, we obtain
\begin{equation}\label{eq:first-component-bound}
 h_1r_1\geqslant a\left(|S|+\sum_{j=2}^c\frac{h_jr_j}{b_j}\right).
\end{equation}
On the other hand, multiplying \eqref{eq:schur-bound} by
$\mu_n-r_1>0$ gives
\[
 |S|(\mu_n-r_1)
 \geqslant h_1r_1
      +\sum_{j=2}^c\frac{h_jr_j(\mu_n-r_1)}{\mu_n-r_j}.
\]
Insert \eqref{eq:first-component-bound} and subtract $a|S|$ from both
sides. Since $\mu_n-r_1-a=\mu_n-\mu_2$, this yields
\begin{equation}\label{eq:combined-cut-bound}
 |S|(\mu_n-\mu_2)
 \geqslant\sum_{j=2}^c h_jr_j
 \left(\frac{a}{b_j}+\frac{\mu_n-r_1}{\mu_n-r_j}\right).
\end{equation}
Using $\mu_n-r_1=(\mu_n-r_j)+a+b_j$, we rewrite each factor in
parentheses as
\[
 \frac{a}{b_j}+\frac{\mu_n-r_1}{\mu_n-r_j}
 =1+\frac{a}{b_j}+\frac{a+b_j}{\mu_n-r_j}
 =\left(1+\frac{a}{b_j}\right)\left(1+\frac{b_j}{\mu_n-r_j}\right)
 =\left(1+\frac{a}{b_j}\right)\frac{\mu_n-\mu_2}{\mu_n-r_j}.
\]
Since $\mu_n-\mu_2>0$ by Lemma~\ref{lem:basic-bounds}, dividing
\eqref{eq:combined-cut-bound} by $\mu_n-\mu_2$ gives
\[
 |S|\geqslant\sum_{j=2}^c\frac{h_jr_j}{\mu_n-r_j}
 \left(1+\frac{a}{b_j}\right).
\]
By \eqref{eq:denominator-bound},
$h_j/(\mu_n-r_j)\geqslant1/(\mu_n-\delta)$. Together with $a\geqslant b_j$,
$r_j>\mu_2$, and $c\geqslant2$, this gives
\[
 |S|(\mu_n-\delta)
 \geqslant\sum_{j=2}^c r_j\left(1+\frac{a}{b_j}\right)
 \geqslant2\sum_{j=2}^c r_j
 >2(c-1)\mu_2
 \geqslant c\mu_2,
\]
which contradicts the assumption. This proves \eqref{eq:strict-cut}.
\end{proof}

The following proposition gives the cut inequality together with the
constraints that equality imposes on every component.

\begin{proposition}\label{prop:cut-bound-equality}
Every vertex cut $S$ of $G$ satisfies
\begin{equation}\label{eq:cut-bound}
 \frac{|S|}{c(G-S)}\geqslant\frac{\mu_2}{\mu_n-\delta}.
\end{equation}
If equality holds, then every component of $G-S$ is a single vertex of
degree $\delta$, and $\mu_2=\delta$.
\end{proposition}

\begin{proof}
If $r_i<\mu_2$ for some $i$, then
Proposition~\ref{prop:exceptional-strictness} gives strict inequality in
\eqref{eq:cut-bound}. Otherwise $r_i\geqslant\mu_2$ for every $i$, and
\eqref{eq:sum-bound} gives
$|S|(\mu_n-\delta)\geqslant\sum_{i=1}^c r_i\geqslant c\mu_2$. This
proves \eqref{eq:cut-bound}.

Suppose that equality holds. Then $r_i\geqslant\mu_2$ for every $i$ by
Proposition~\ref{prop:exceptional-strictness}, and
\eqref{eq:schur-bound} and \eqref{eq:denominator-bound} give
\[
 c\mu_2=|S|(\mu_n-\delta)
 \geqslant\sum_{i=1}^c
          \frac{(\mu_n-\delta)h_ir_i}{\mu_n-r_i}
 \geqslant\sum_{i=1}^c r_i
 \geqslant c\mu_2.
\]
Hence all these inequalities are equalities. The last one gives
$r_i=\mu_2$ for every $i$. Since $r_i>0$, the middle one gives
$\mu_n-r_i=(\mu_n-\delta)h_i$, so both inequalities in
\eqref{eq:denominator-bound} are equalities. By
Lemma~\ref{lem:upper-cut}, this forces $h_i=1$, and hence
$r_i=\delta$, because $\delta-h_i + 1 \le r_i = \mu_2 \le \delta$. For a singleton component, $r_i$ is the degree of its
vertex. Hence every component of $G-S$ is a single vertex of degree
$\delta$, and $\mu_2=r_i=\delta$.
\end{proof}

Taking the minimum of \eqref{eq:cut-bound} over all vertex cuts $S$
recovers \eqref{eq:haemers-bound}. Moreover,
Proposition~\ref{prop:cut-bound-equality} reduces the equality problem
to vertex cuts $S$ such that $V(G)\setminus S$ is an independent set of
vertices of degree $\mu_2$.

\section{The equality characterization}\label{sec:characterization}

We first show that, in an equality graph, the vertices outside a
toughness cut form an independent set joined to every vertex of the cut.
The minimality of the toughness cut forces the common neighborhood given
by Lemma~\ref{lem:common-neighborhood} to be the whole cut.

\begin{lemma}\label{lem:equality-cut-structure}
Let $G$ be an equality graph of order $n$, and let $S$ be a toughness
cut of $G$. Then $V(G)\setminus S$ is a maximum independent set, every
vertex of $V(G)\setminus S$ is adjacent to every vertex of $S$, and
\begin{equation}\label{eq:cut-structure-parameters}
 |S|=\delta=\mu_2,
 \qquad
 |V(G)\setminus S|=n-\delta,
 \qquad
 \mu_n=n.
\end{equation}
\end{lemma}

\begin{proof}
Put $U=V(G)\setminus S$ and $m=|U|$. Since $S$ is a toughness cut,
$|S|/c(G-S)=t(G)=\mu_2/(\mu_n-\delta)$. By
Proposition~\ref{prop:cut-bound-equality}, every component of $G-S$ is a
single vertex of degree $\delta$, and $\mu_2=\delta$. Thus $U$ is an
independent set of vertices of degree $\delta$, $c(G-S)=m\geqslant2$,
and $t(G)=|S|/m$. By Lemma~\ref{lem:common-neighborhood}, the vertices
of $U$ have a common neighborhood $W\subseteq S$, and $|W|=\delta$.

Suppose that $S\setminus W$ is nonempty. Deleting $W$ isolates the $m$
vertices of $U$, and the vertices of $S\setminus W$ form at least one
further component. Hence
\[
 t(G)\leqslant\frac{|W|}{c(G-W)}
 \leqslant\frac{\delta}{m+1}
 \leqslant\frac{|S|}{m+1}
 <\frac{|S|}{m}=t(G),
\]
a contradiction. Therefore $W=S$. Thus every vertex of $U$ is adjacent
to every vertex of $S$, $|S|=\delta=\mu_2$, and $m=n-\delta$.

The equality $t(G)=\mu_2/(\mu_n-\delta)$ now reads
$\delta/m=\delta/(\mu_n-\delta)$. Since $\delta>0$, we obtain
$\mu_n=\delta+m=n$. Finally, \eqref{eq:alpha-chain} gives
$\alpha(G)\leqslant n-\delta=m$, so the independent set $U$ is maximum.
\end{proof}

\begin{proof}[Proof of Theorem~\ref{thm:main}]
\emph{\textup{(i)} implies \textup{(ii)} and \textup{(iii)}.}
Assume \textup{(i)}, let $S$ be a toughness cut, and put
$U=V(G)\setminus S$ and $m=n-\delta$. By
Lemma~\ref{lem:equality-cut-structure}, $U$ is a maximum independent set
of size $m\geqslant2$, and every vertex of $U$ is adjacent to every
vertex of $S$. Hence $\overline G[U]\cong K_m$ is a component of
$\overline G$, and
\[
 \overline G\cong K_m\sqcup F,
 \qquad F=\overline G[S],
\]
where $F$ has $|S|=\delta\geqslant1$ vertices. By
Lemma~\ref{lem:complement} and $\mu_2=\delta=n-m$, we have
$\rhL(\overline G)=n-\mu_2=m$. Since $L(\overline G)$ is block diagonal
with blocks $L(K_m)$ and $L(F)$, it follows that $\rhL(F)\leqslant m$,
which proves \textup{(iii)}. Moreover, $\alpha(G)=|U|=m=n-\mu_2$, which
proves \textup{(ii)}.

\emph{\textup{(iii)} implies \textup{(i)}.}
Assume \textup{(iii)}. Since $\rhL(K_m)=m$ and $\rhL(F)\leqslant m$, we
have $\rhL(\overline G)=m$. As $\overline G$ is disconnected,
Lemma~\ref{lem:complement} gives
\[
 \mu_2=n-m
 \qquad\text{and}\qquad
 \mu_n=n.
\]
The $m$ vertices corresponding to $K_m$ have degree $n-m$ in $G$, so
$\delta\leqslant n-m=\mu_2\leqslant\delta$, and hence $\delta=n-m$.
Deleting the other $n-m$ vertices leaves $m\geqslant2$ isolated vertices,
so
\[
 t(G)\leqslant\frac{n-m}{m}
       =\frac{\mu_2}{\mu_n-\delta}.
\]
The reverse inequality is Proposition~\ref{prop:cut-bound-equality}.
This proves \textup{(i)} and \eqref{eq:equality-parameters}. Finally,
$\alpha(G)\leqslant n-\delta=m$ by \eqref{eq:alpha-chain}, and the $m$
vertices corresponding to $K_m$ form an independent set in $G$; hence
$m=\alpha(G)$.

\emph{\textup{(ii)} implies \textup{(iii)}.}
Assume \textup{(ii)}. By \eqref{eq:alpha-chain},
$\mu_2=\delta=n-\alpha(G)$. Let $U$ be a maximum independent set, and
put $m=\alpha(G)$. Every vertex of $U$ has all its neighbors outside
$U$, so its degree is at most $n-m=\delta$. Hence it has degree exactly
$\delta$ and is adjacent to every vertex outside $U$. Therefore
\[
 \overline G\cong K_m\sqcup F,
 \qquad F=\overline G[V(G)\setminus U].
\]
Here $m\geqslant2$ because $G$ is noncomplete, and $F$ has at least one
vertex because $G$ is connected. As before, Lemma~\ref{lem:complement}
gives $\rhL(\overline G)=n-\mu_2=m$, and hence $\rhL(F)\leqslant m$.
This proves \textup{(iii)} and completes the proof.
\end{proof}

We now compare the equality cases of the three Laplacian bounds.

\begin{corollary}\label{cor:three-bounds}
For a connected noncomplete graph $G$ of order $n$, the following three
equalities are equivalent:
\[
 t(G)=\frac{\mu_2}{\mu_n-\delta},
 \qquad
 t(G)=\frac{\mu_n\mu_2}{n(\mu_n-\delta)},
 \qquad
 t(G)=\frac{\mu_2}{\mu_n-\mu_2}.
\]
In particular, the bounds \eqref{eq:gu-haemers-first},
\eqref{eq:gu-haemers-second}, and \eqref{eq:haemers-bound} have the same
equality graphs, namely those described in Theorem~\ref{thm:main}.
\end{corollary}

\begin{proof}
By Lemma~\ref{lem:basic-bounds} and \eqref{eq:haemers-bound},
\[
 \frac{\mu_n\mu_2}{n(\mu_n-\delta)}
 \leqslant\frac{\mu_2}{\mu_n-\delta}
 \leqslant t(G)
 \qquad\text{and}\qquad
 \frac{\mu_2}{\mu_n-\mu_2}
 \leqslant\frac{\mu_2}{\mu_n-\delta}
 \leqslant t(G).
\]
Hence each of the last two equalities implies the first. Conversely, the
first equality gives $\mu_n=n$ and $\mu_2=\delta$ by
Theorem~\ref{thm:main}, and then the three ratios coincide.
\end{proof}

\begin{remark}\label{rem:numerical-conditions}
By Theorem~\ref{thm:main}, inequality \eqref{eq:haemers-bound} is strict
whenever $\mu_2<\delta$ or $\mu_n<n$; by Lemma~\ref{lem:complement}, the
latter holds exactly when $\overline G$ is connected. The identities
$\mu_2=\delta$ and $\mu_n=n$ are necessary for equality but not
sufficient: the graph in Example~\ref{ex:necessary-not-sufficient}
satisfies both identities, and \eqref{eq:haemers-bound} is strict for
it. The complete component in
Theorem~\ref{thm:main}\textup{(iii)} is therefore an essential part of
the characterization.
\end{remark}

\begin{example}\label{ex:necessary-not-sufficient}
Let $G=\overline{P_3\sqcup P_3}$, where $P_3$ is the path on three
vertices. Since $L(P_3)$ has eigenvalues $0,1,3$,
Lemma~\ref{lem:complement} gives $\mu_2(G)=6-3=3$ and $\mu_6(G)=6$.
Moreover, $\delta(G)=3$, so the right-hand side of
\eqref{eq:haemers-bound} equals $1$. We show that $t(G)=3/2$.

The graph $G$ is the join $(K_2\sqcup K_1)\vee(K_2\sqcup K_1)$. If $G-S$
contains a vertex from each join factor, then $G-S$ is connected, since
these two vertices are adjacent and every other vertex of $G-S$ is
adjacent to one of them. Hence every vertex cut contains one of the two
factors and has at least three vertices. Moreover, choosing one vertex
from each component of $G-S$ gives an independent set, so
$c(G-S)\leqslant\alpha(G)=2$. Thus every vertex cut has ratio at least
$3/2$, and deleting either join factor attains this value. Accordingly,
$\overline G=P_3\sqcup P_3$ has no complete component of order at least
two, and condition \textup{(iii)} of Theorem~\ref{thm:main} fails.
\end{example}

\section{Regular graphs and toughness cuts}\label{sec:consequences}

We now derive Corollary~\ref{cor:regular-introduction} and determine all
toughness cuts of an equality graph. Under regularity, every component of
the complement is a copy of $K_m$; in general, the components of the
complement isomorphic to $K_m$ correspond exactly to the toughness cuts.

\begin{proof}[Proof of Corollary~\ref{cor:regular-introduction}]
Let $G$ be a regular equality graph of order $n$. By
Theorem~\ref{thm:main}, there is an integer $m\geqslant2$ such that
$\delta(G)=n-m$ and $\rhL(\overline G)=n-\mu_2(G)=m$. Since $G$ is
regular of degree $n-m$, its complement is $(m-1)$-regular. Let $C$ be a
connected component of $\overline G$. If $C$ were noncomplete,
Lemma~\ref{lem:basic-bounds} applied to $C$ would give
\[
 \rhL(C)>\delta(C)+1=m,
\]
contrary to $\rhL(\overline G)=m$. Hence every component of
$\overline G$ is complete and $(m-1)$-regular, that is, isomorphic to
$K_m$. Since $G$ is connected, $\overline G$ has at least two
components.

Conversely, if $\overline G$ is the disjoint union of $r\geqslant2$
copies of $K_m$ with $m\geqslant2$, then $\overline G$ satisfies
Theorem~\ref{thm:main}\textup{(iii)}, so $G=K_{m,\ldots,m}$ is an
equality graph. Its order is $rm$, and \eqref{eq:equality-parameters}
gives $t(G)=(rm-m)/m=r-1$.
\end{proof}

\begin{theorem}\label{thm:all-cuts}
Let $G$ be an equality graph of order $n$, and put $m=\alpha(G)$. For a
subset $S\subseteq V(G)$, the following statements are equivalent.
\begin{enumerate}
 \item[\textup{(i)}] $S$ is a toughness cut of $G$.
 \item[\textup{(ii)}] $V(G)\setminus S$ is a maximum independent set of
 $G$.
 \item[\textup{(iii)}] $V(G)\setminus S$ is the vertex set of a
 component of $\overline G$ isomorphic to $K_m$.
\end{enumerate}
Consequently, the number of toughness cuts of $G$ equals the number of
components of $\overline G$ isomorphic to $K_m$. Every toughness cut has
size $n-m=\delta(G)$ and leaves exactly $m$ components.
\end{theorem}

\begin{proof}
The implication \textup{(i)}$\Rightarrow$\textup{(ii)} is part of
Lemma~\ref{lem:equality-cut-structure}.

For \textup{(ii)}$\Rightarrow$\textup{(iii)}, let $U=V(G)\setminus S$ be
a maximum independent set, so that $|U|=m$. By Theorem~\ref{thm:main},
$\delta(G)=n-m$. Every vertex of $U$ has all its neighbors in $S$ and at
least $n-m=|S|$ of them, so it is adjacent to every vertex of $S$. In
$\overline G$, the set $U$ therefore induces $K_m$ and sends no edges to
$S$; hence it is the vertex set of a component.

For \textup{(iii)}$\Rightarrow$\textup{(i)}, the graph $G-S$ consists of
$m\geqslant2$ isolated vertices, and $|S|=n-m$. Its ratio is
$(n-m)/m=t(G)$ by \eqref{eq:equality-parameters}, so $S$ is a toughness
cut.

Distinct components of $\overline G$ have distinct complements in
$V(G)$, which gives the counting statement. The remaining assertions
follow from \textup{(iii)} and \eqref{eq:equality-parameters}.
\end{proof}

In particular, the maximum independent sets of an equality graph are
pairwise disjoint, there are at most $\lfloor n/m\rfloor$ of them, and
the toughness cut is unique if and only if $\overline G$ has exactly one
component isomorphic to $K_m$.

\begin{corollary}\label{cor:minimum-cuts}
Every equality graph $G$ satisfies $\kappa(G)=\delta(G)=\mu_2(G)$, and
every toughness cut of $G$ is a minimum vertex cut.
\end{corollary}

\begin{proof}
By Lemma~\ref{lem:basic-bounds} and Theorem~\ref{thm:main}, every vertex
cut has size at least $\mu_2(G)=\delta(G)$. By
Theorem~\ref{thm:all-cuts}, every toughness cut has size $\delta(G)$.
\end{proof}

\begin{example}\label{ex:minimum-not-toughness}
A minimum vertex cut of an equality graph need not be a toughness cut.
Let $G=\overline K_3\vee(K_2\sqcup K_1)$. Its complement is
$K_3\sqcup P_3$, and $\rhL(P_3)=3$, so Theorem~\ref{thm:main} applies
with $m=3$ and gives $\delta(G)=3$ and $t(G)=1$. Deleting the three
vertices of $K_2\sqcup K_1$ leaves three isolated vertices; this is a
toughness cut. Deleting instead the three vertices of $\overline K_3$
leaves $K_2\sqcup K_1$, which has two components. This minimum vertex
cut has ratio $3/2$.
\end{example}

\section{Examples and recognition}\label{sec:examples}

\subsection{The join description}\label{subsec:join}

By Theorem~\ref{thm:main}, applied with $H=\overline F$, a connected
noncomplete graph $G$ is an equality graph if and only if
\begin{equation}\label{eq:join-description}
 G\cong\overline K_m\vee H,
 \qquad m\geqslant2,
 \qquad |V(H)|\geqslant1,
 \qquad \rhL(\overline H)\leqslant m.
\end{equation}
In this case, $m=\alpha(G)$ and $|V(H)|=\delta(G)$. If $H$ has order
$s\geqslant2$, then Lemma~\ref{lem:complement} rewrites the condition on
$H$ as
\begin{equation}\label{eq:join-gap}
 \mu_2(H)\geqslant s-m;
\end{equation}
this includes disconnected graphs $H$, for which $\mu_2(H)=0$. If $s=1$,
the condition in \eqref{eq:join-description} holds trivially. With
$s=\delta$ and $m=n-\delta$, condition \eqref{eq:join-gap} is the form
in which Huang, Das, and Zhu~\cite{HuangDasZhu2022} describe their
equality family.

The condition $\rhL(\overline H)\leqslant m$ holds automatically when
$|V(H)|\leqslant m$, because the Laplacian spectral radius of a graph is
at most its order (Lemma~\ref{lem:complement}). Hence every graph on at
most $m$ vertices occurs as $H$. Adding edges to $H$ also preserves
equality.

\begin{proposition}\label{prop:edge-addition}
Let $G=\overline K_m\vee H$ satisfy \eqref{eq:join-description}, and let
$H'$ be obtained from $H$ by adding edges on the same vertex set. Then
$\overline K_m\vee H'$ is an equality graph with the same toughness and
independence number as $G$.
\end{proposition}

\begin{proof}
The graph $\overline{H'}$ is obtained from $\overline H$ by deleting
edges. Deleting an edge $uv$ subtracts the positive semidefinite matrix
$(\ee_u-\ee_v)(\ee_u-\ee_v)^{\mathsf T}$ from the Laplacian matrix, so
the maximum Rayleigh quotient gives
$\rhL(\overline{H'})\leqslant\rhL(\overline H)\leqslant m$. By
Theorem~\ref{thm:main}, both graphs are equality graphs with
independence number $m$ and toughness $|V(H)|/m$.
\end{proof}

\begin{example}\label{ex:multipartite}
Let $G=K_{n_1,\ldots,n_r}$ with $r\geqslant2$ and
$m=\max_in_i\geqslant2$. Its complement is
$K_{n_1}\sqcup\cdots\sqcup K_{n_r}$. Removing one component isomorphic
to $K_m$ leaves a disjoint union of complete graphs of order at most $m$,
whose Laplacian spectral radius is at most $m$. Hence $G$ is an equality
graph with $t(G)=(n-m)/m$. By Theorem~\ref{thm:all-cuts}, its toughness
cuts are exactly the complements of its parts of size $m$.
\end{example}

\begin{example}\label{ex:split}
For $m\geqslant2$ and $s\geqslant1$, the graph $\overline K_m\vee K_s$
is an equality graph with toughness $s/m$, since its complement is
$K_m\sqcup\overline K_s$ and $\rhL(\overline K_s)=0$. In particular,
every positive rational number is the toughness of an equality graph:
for positive integers $a$ and $b$, the graph
$\overline K_{2b}\vee K_{2a}$ has toughness $a/b$.
\end{example}

In \eqref{eq:join-description}, the factor $\overline K_m$ is a maximum
independent set of $G$. For example, $K_{2,3}=\overline K_3\vee\overline
K_2$ satisfies \eqref{eq:join-description} with $m=3$, because
$\rhL(K_2)=2\leqslant3$; the decomposition
$K_{2,3}=\overline K_2\vee\overline K_3$ does not, because
$\rhL(K_3)=3>2$.

\subsection{Noncomplete components of the complement}\label{subsec:noncomplete}

The following equality graphs are not regular, and their complements have
noncomplete components of order larger than the independence number. For
$q\geqslant1$, let $P_q$ be the path on $q$ vertices, and for
$q\geqslant3$, let $C_q$ be the cycle on $q$ vertices.

\begin{example}\label{ex:cycle-five}
The graph $G=\overline{K_4\sqcup C_5}$ is an equality graph. Indeed, the
characteristic polynomial of $L(C_5)$ is $x(x^2-5x+5)^2$, so
\[
 \rhL(C_5)=\frac{5+\sqrt5}{2}<4.
\]
Theorem~\ref{thm:main} gives $n=9$, $\mu_2=\delta=5$, $\mu_n=9$, and
$t(G)=5/4$. The four vertices corresponding to $K_4$ have degree five in
$G$, and the other five vertices have degree six. Thus $G$ is not
regular, and the component $C_5$ of $\overline G$ has order larger than
$\alpha(G)=4$.
\end{example}

\begin{example}\label{ex:paths}
For every $q\geqslant1$, the graph $G_q=\overline{K_4\sqcup P_q}$ is an
equality graph with $t(G_q)=q/4$ and $\alpha(G_q)=4$. For $q=1$, the
graph $P_1$ has Laplacian spectral radius zero; for $q\geqslant2$,
Lemma~\ref{lem:degree-radius} gives $\rhL(P_q)\leqslant4$. The same
holds with $P_q$ replaced by any nonempty disjoint union of paths and
cycles. Hence the orders of the noncomplete components of the complement
of an equality graph are not bounded in terms of its independence number.
\end{example}

For independence number two or three, the condition in
Theorem~\ref{thm:main}\textup{(iii)} reduces to a finite list of
admissible components.

\begin{proposition}\label{prop:small-independence}
Let $G$ be a connected noncomplete graph.
\begin{enumerate}
 \item[\textup{(i)}] $G$ is an equality graph with $\alpha(G)=2$ if and
 only if every component of $\overline G$ is $K_1$ or $K_2$ and at least
 one component is $K_2$.
 \item[\textup{(ii)}] $G$ is an equality graph with $\alpha(G)=3$ if and
 only if every component of $\overline G$ is one of
 $K_1,K_2,P_3,K_3$ and at least one component is $K_3$.
\end{enumerate}
\end{proposition}

\begin{proof}
\textup{(i)} By Theorem~\ref{thm:main}, $G$ is an equality graph with
$\alpha(G)=2$ if and only if $\overline G$ has a component $K_2$ and the
remaining components have Laplacian spectral radius at most two. By
Lemma~\ref{lem:degree-radius}, a graph with Laplacian spectral radius at
most two has maximum degree at most one, so its components are $K_1$ and
$K_2$. Conversely, $\rhL(K_1)=0$ and $\rhL(K_2)=2$.

\textup{(ii)} Similarly, $G$ is an equality graph with $\alpha(G)=3$ if
and only if $\overline G$ has a component $K_3$ and the remaining
components have Laplacian spectral radius at most three. By
Lemma~\ref{lem:degree-radius}, such components have maximum degree at
most two, so each is a path or a cycle. None of them contains $P_4$ as a
subgraph: the vector $(1,-2,2,-1)^{\mathsf T}$ on the consecutive
vertices of $P_4$, extended by zeros, has Rayleigh quotient at least
$34/10>3$ in any graph containing this path, because additional edges
contribute nonnegative terms to \eqref{eq:laplacian-energy}. Hence the
admissible components are exactly $K_1$, $K_2$, $P_3$, and $K_3$. Each
of them has at most three vertices, so its Laplacian spectral radius is
at most three by Lemma~\ref{lem:complement}, with $\rhL(K_1)=0$.
\end{proof}

\subsection{A recognition criterion}\label{subsec:recognition}

By \eqref{eq:equality-parameters}, the parameter $m$ in
Theorem~\ref{thm:main} is determined by the minimum degree. This gives a
test for equality that requires neither the independence number nor a
toughness cut.

\begin{proposition}\label{prop:recognition}
Let $G$ be a connected noncomplete graph of order $n$, and put
$m=n-\delta(G)$. Then $G$ is an equality graph if and only if
$\overline G$ has a component $C\cong K_m$ such that the matrix
$mI-L(F)$ is positive semidefinite, where $F=\overline G-V(C)$. In this
case, the toughness cuts of $G$ are exactly the complements of the
components of $\overline G$ isomorphic to $K_m$.
\end{proposition}

\begin{proof}
By \eqref{eq:equality-parameters}, every integer $m$ satisfying
Theorem~\ref{thm:main}\textup{(iii)} equals $n-\delta(G)$. Hence $G$ is
an equality graph if and only if $\overline G$ has a component
$C\cong K_m$ with $m=n-\delta(G)$ such that $F=\overline G-V(C)$
satisfies $\rhL(F)\leqslant m$; the graph $F$ has at least one vertex
because $G$ is connected. Since $L(F)$ is symmetric,
$\rhL(F)\leqslant m$ if and only if $mI-L(F)$ is positive semidefinite.
The statement about toughness cuts is Theorem~\ref{thm:all-cuts}.
\end{proof}

\end{document}